\documentclass[12pt]{amsart}
\usepackage{tikz}
\usetikzlibrary{matrix,decorations.pathreplacing}

\usepackage{amsmath,amsthm,amssymb}
\font\teneufm=eufm10
\font\seveneufm=eufm7
\font\fiveeufm=eufm5
\newfam\frakturfam

\textfont\frakturfam=\teneufm
\scriptfont\frakturfam=\seveneufm
\scriptscriptfont\frakturfam=\fiveeufm

\theoremstyle{definition}
\newtheorem{definition}{Definition}[section]
\newtheorem{remark}[definition]{Remark}

\theoremstyle{plain}
\newtheorem{lemma}[definition]{Lemma}
\newtheorem{prp}[definition]{Proposition}
\newtheorem{corollary}[definition]{Corollary}
\newtheorem{theorem}[definition]{Theorem}

\usepackage{xcolor}
\usepackage{array}

\begin{document}

\title{Polynomial identities of the two-dimensional simple right-symmetric algebra}

\author{Nurlan Ismailov}

\address{Astana IT University,
Mangilik El avenue, 55/11, Business center EXPO, block C1,
Astana, 010000, Kazakhstan and Institute of Mathematics and Mathematical Modeling, Almaty, 050010, Kazakhstan}

\email{nurlan.ismail@gmail.com}

\author{Ualbai Umirbaev}

\address{Department of Mathematics,
 Wayne State University, Detroit, MI 48202, USA; Center for Theoretical and Applied Algebra, Al-Farabi Kazakh National University, Almaty, 050040
and Institute of Mathematics and Mathematical Modeling, Almaty, 050010, Kazakhstan}

\email{umirbaev@wayne.edu}

\subjclass[2020]{Primary 17D25, 17A30. Secondary 17A50}

\keywords{right-symmetric algebras, pre-Lie algebras, two-dimensional simple algebras, polynomial identities, $T$-ideal}

\maketitle

\begin{abstract} 
We describe a basis for the polynomial identities satisfied by the unique two-dimensional simple right-symmetric algebra $B$ over a field of characteristic zero. Modulo the right-symmetric identity, this basis consists of two irreducible identities of degree three.  We construct linear bases for the multilinear components of the relatively free algebra of the variety 
$\mathcal V=\operatorname{var}(B)$ and determine their decompositions into Specht modules.  Using the resulting $S_n$-module decompositions, we prove that the commutator subspace is right nilpotent in every proper subvariety of $\mathcal V$. Moreover, in any proper subvariety satisfying a nontrivial identity of degree $n\ge2$, the index of right nilpotence of the commutator subspace is at most $2n$.

\end{abstract}

\maketitle
\section{Introduction}

  Let $\mathbb{F}\{X\}$ be the free nonassociative algebra over a field $\mathbb{F}$
on the countable set $X=\{x_1,x_2,\ldots\}$.
A polynomial $f(x_1,\ldots,x_n)\in \mathbb{F}\{X\}$ is said to be an
\emph{identity} of an $\mathbb{F}$-algebra $A$ if
$$
f(a_1,\ldots,a_n)=0
\qquad\text{for all }a_1,\ldots,a_n\in A.
$$
An ideal $I$ of $\mathbb{F}\{X\}$ is called a \emph{$T$-ideal} if it
is invariant under every endomorphism of $\mathbb{F}\{X\}$. The set of all polynomial identities of $A$, denoted by
$\operatorname{Id}(A)$, forms a $T$-ideal of $\mathbb{F}\{X\}$. Any set that generates $\operatorname{Id}(A)$ as a $T$-ideal is called a {\em basis of identities} of $A$. 

The study of identities of simple algebras, particularly those of “nontrivial” simple algebras of minimal dimension in each variety, is of great importance in ring theory. Shestakov and Zaicev  \cite{ShestakovZaicev2011} proved that, over an algebraically closed field,  two simple  finite-dimensional algebras satisfy the same polynomial identities  if and only if  they are isomorphic. The algebra $M_2$ of $2\times2$ matrices has the smallest dimension among noncommutative simple associative algebras, while the three-dimensional Lie algebra $\mathfrak{sl}_2$
of traceless $2\times2$ matrices has the smallest dimension among simple Lie algebras. Razmyslov \cite{Razmyslov1973} proved that, over a field of characteristic zero, both algebras have finite bases of polynomial identities.  He also established
the Specht property for the variety generated by
$\mathfrak{sl}_2$ \cite{Razmyslov1974}. Filippov \cite{Filippov1981} showed that all identities of $\mathfrak{sl}_2$ can be reduced to a single
identity. Drensky \cite{Drensky1981} subsequently obtained a minimal basis of
identities for $M_2$ and another explicit basis for $\mathfrak{sl}_2$,
using the representation theory of symmetric and general linear
groups.

Grassmann (exterior) algebras play an important role in the theory of polynomial identities. A basis for their identities was studied in \cite{KrakowskiRegev1973,Latyshev1963}. 

An algebra $A$ is called {\em right-symmetric} (pre-Lie) if it satisfies
\begin{equation}\label{rsym}(a,b,c)=(a,c,b), 
\end{equation}
where $(a,b,c)=(a b) c-a (b c)$ is the associator.  Right-symmetric algebras arise in many different areas of mathematics and physics \cite{Burde2006}. 
Every right-symmetric
algebra is Lie admissible: the commutator $[a,b]=ab-ba$ defines a Lie
algebra structure on $A$.

In \cite{IsmailovUmirbaev2024}, the authors of the present paper  
constructed a finite-dimensional metabelian right-symmetric algebra
over an arbitrary field whose identities do not have a finite basis.
Thus, the variety of right symmetric algebras does not have the Specht property.

 A right-symmetric algebra satisfying the additional identity 
$$
a (b c)=b (a c)
$$ is called a Novikov algebra. The class of Novikov algebras is an important subclass of right-symmetric algebras. The polynomial identities of Novikov algebras have been extensively studied.   Dotsenko, Ismailov, and Umirbaev \cite{DIU2023} proved that (a) every Novikov algebra satisfying a nontrivial polynomial identity over a field of characteristic zero is right-associator nilpotent and (b) the variety of Novikov algebras over a field of characteristic zero has the Specht property. Zelmanov \cite{Zelmanov1987} proved that every finite-dimensional simple Novikov algebra over an algebraically closed field of characteristic zero has dimension one.  In \cite{SantosKuzminLopatin2025}, Ferreira dos Santos, Kuz’min, and Lopatin studied minimal bases of polynomial identities and bases of the corresponding relatively free algebras for all two-dimensional complex Novikov algebras.  Giambruno, Mishchenko, and Zaicev~\cite{GMZ2007} proved that, over a field of characteristic zero, the codimension sequence of a two-dimensional algebra is either bounded by $n+1$ or has PI-exponent $2$. Giambruno, Shestakov, and Zaicev~\cite{GSZ2011} studied codimension growth for broader classes of finite-dimensional nonassociative algebras.

A powerful approach to polynomial identities is to study the multilinear components of relatively free algebras and their decompositions into irreducible modules over symmetric groups. The dimensions of these components form the codimension sequence, which measures the growth of the relatively free algebra in multilinear degree.

In this paper, we study the polynomial identities of the noncommutative simple right-symmetric algebra of minimal dimension over a field of characteristic zero.    
In \cite{Burde1998}, Burde classified simple right-symmetric algebras in low dimensions over the complex numbers. Up to isomorphism, there is only one two-dimensional simple right-symmetric algebra $B=\langle e_1,e_2\rangle$ over $\mathbb{C}$, with multiplication given by
$$
e_1e_1=2e_1,\qquad e_1e_2=0,\qquad
e_2e_1=e_2,\qquad e_2e_2=e_1
$$
(see \cite[Proposition~4.1]{Burde1998}, using the opposite convention).

Pozhidaev, Umirbaev, and
Zhelyabin studied simple left-symmetric algebras in \cite{PUZ2022}. The opposite algebra $B^{\mathrm{op}}$ is left-symmetric and satisfies
$$[[x,y],[z,t]]=0,\qquad [x,y]([z,t]u)=0.$$
Thus, it belongs to the class of simple left-symmetric algebras $\mathcal{A}$ over an algebraically closed field studied in \cite{PUZ2022}.

We prove that all polynomial identities of $B$ follow from the
right-symmetric identity and the following two irreducible identities of degree three:
$$
 a[b,c]+b[c,a]+c[a,b]=0
$$
and
$$
 3a[a,b]-2[a,b]a-b(aa)+a(ab)=0.
$$
Consequently, these three identities define the variety
$\mathcal V=\operatorname{var}(B)$. We construct a linear basis of
each multilinear component $P_n(\mathcal V)$ and decompose it
into irreducible $S_n$-modules. Specifically, for
$n\geq2$,
$$
 P_n(\mathcal V)\cong
 (n-1)S^{(n)}\oplus
 \bigoplus_{j=1}^{\lfloor n/2\rfloor}
 (n-2j+1)S^{(n-j,j)},
$$
where $S^\lambda$ is the Specht module (irreducible $S_n$-module) corresponding to the
partition $\lambda$ of $n$. In particular,
$$
 c_1(\mathcal V)=1,\qquad
 c_n(\mathcal V)=2^n-2\quad(n\ge2),$$
and the PI-exponent of $\mathcal V$ equals $2$.

The cocharacter and codimension sequences of $\mathcal V$ coincide with
those of the variety of all bicommutative algebras \cite {Dzh-Ism-Tul2011} (see also \cite{Drensky2019}). An algebra satisfying the identities 
$$
 (ab)c=(ac)b,\qquad a(bc)=b(ac)
$$
is called bicommutative. Drensky \cite{Drensky2019}  classified two-dimensional bicommutative algebras and showed that one of them generates the entire variety of bicommutative algebras. He also proved that the codimensions of every proper subvariety of the variety of bicommutative algebras are polynomially bounded.

We prove that every proper subvariety $\mathcal W\subsetneq\mathcal V$
satisfies an identity of the form 
$$
 [x_1,y_1]\cdots[x_s,y_s]=0
$$
for some $s\geq1$, where the product is right-normed. The integer $s$
depends only on $\mathcal W$. More precisely, if $\mathcal W$ satisfies an identity of degree $n\ge2$ that does not hold in $\mathcal V$, then $s\leq2n$.
For an identity that generates an irreducible module corresponding to a partition $(\lambda_1,\lambda_2)$ with $\lambda_2>0$, we get the sharper bound
$s\le2\lambda_1+\lambda_2$. The proof proceeds by relating irreducible multilinear identities first to powers of commutators and then to right-normed products of commutators.

The paper is organized as follows.  In Section 2,  we establish the polynomial and operator identities used throughout
the paper. In Section 3, we construct bases for the multilinear components and obtain a finite basis of identities for $B$. In Section 4, we determine the $S_n$-module structure of each multilinear component. Section 5 is devoted to polynomial identities of proper subvarieties and right-normed products of commutators.

\section{Polynomial identities of $B$}
Let $\mathbb{F}$ be a field of characteristic zero. The right-symmetric identity \eqref{rsym} is equivalent to the identity
\begin{equation}\label{rsym-commutator} 
a[b,c]=(ab)c-(ac)b
\end{equation}
where $[a,b]=ab-ba$.
Recall that every right-symmetric algebra is Lie admissible: 
$$Jac(a,b,c)=[[a,b],c]+[[b,c],a]+[[c,a],b]=0.$$
Let $B=span_{\mathbb{F}}\{e_1,e_2\}$ be the
two-dimensional simple right-symmetric algebra over $\mathbb{F}$
with the following multiplication table \cite[Proposition 4.1]{Burde1998}:
$$e_1e_1=2e_1, \qquad e_1e_2=0, \qquad
e_2e_1=e_2, \qquad e_2e_2=e_1.
$$
Note that $[B,B]=span_{\mathbb{F}}\{e_2\}$.
\begin{prp}\label{B-identities}
  The algebra $B$ satisfies the following identities:
\begin{equation}\label{half-lie-1}
    a[b,c]+b[c,a]+c[a,b]=0,
\end{equation}
\begin{equation}\label{irred-(2,1)}
    3a[a,b]-2[a,b]a-b(aa)+a(ab)=0.
\end{equation} 
\end{prp}

\begin{proof} The expression in \eqref{half-lie-1} is skew-symmetric and trilinear, so it vanishes on the two-dimensional algebra $B$. 
Let
$$
g(x,y,z)=3x[y,z]-2[y,z]x-z(xy)+y(xz).
$$
Since $g$ is trilinear and skew-symmetric in $y$ and $z$, and $\dim B=2$, it suffices to check the following two substitutions. By the multiplication table,
$$
g(e_1,e_1,e_2)=0 \quad\mbox{and}\quad  g(e_2,e_1,e_2)=0.$$
Consequently, for all $a,b\in B$,
$$
3a[a,b]-2[a,b]a-b(aa)+a(ab)=g(a,a,b)=0.
$$\end{proof}

Let $\mathcal{V}$ be  the variety of right-symmetric algebras over $\mathbb{F}$ defined by identities \eqref{half-lie-1} and \eqref{irred-(2,1)}. Below, we give a list of important polynomial identities satisfied by algebras in $\mathcal{V}$.

\begin{lemma} Algebras in $\mathcal{V}$ satisfy the following
identities of degree three:
\begin{equation}\label{degree-3 identity}
3a[b,c]-2[b,c]a-c(ab)+b(ac)=0,
\end{equation}
\begin{equation}\label{mainident1}
4a[b,c]-2[b,c]a-c(ba)+b(ca)=0.
\end{equation}
\end{lemma}

\begin{proof}
Denote the left-hand side of \eqref{degree-3 identity} by $g(a,b,c)$.
The linearization of \eqref{irred-(2,1)} in its repeated variable is
\begin{align*}
 p(a,b,c)={}&3a[b,c]+3b[a,c]-2[a,c]b-2[b,c]a\\
 &+a(bc)+b(ac)-c(ab)-c(ba).
\end{align*}
Expansion gives
$$3g(a,b,c)=p(a,b,c)-p(a,c,b)-2Jac(a,b,c).$$
This proves \eqref{degree-3 identity}.
To obtain \eqref{mainident1}, use \eqref{half-lie-1} in \eqref{degree-3 identity}
$$3a[b,c]-2[b,c]a+b[a,c]+b(ca)-c(ab)=4a[b,c]-2[b,c]a-c(ba)+b(ca).$$
\end{proof}
Write
$$L_a(v)=av,\qquad R_a(v)=va,$$
where operator products act from right to left.
Then identities \eqref{rsym} and \eqref{degree-3 identity} give
\begin{equation}\label{spanning-multiplication-operators}
\begin{aligned}
 L_{ab}&=R_bL_a+L_aL_b-L_aR_b,\\
 R_{ab}&=3L_aL_b-3L_aR_b-2R_aL_b+2R_aR_b+L_bL_a.
\end{aligned}
\end{equation}
\begin{lemma}\label{operator}
If $L=L_x$ and $R=R_x$, then the following operator identity holds:
\begin{equation}\label{spanning-operator-relation}
 2R^2L=2RL^2+RLR-L^2R-2L^3+2LRL.
\end{equation}
\end{lemma}

\begin{proof}
Put $P=L_{x^2}$ and $Q=R_{x^2}$. Identities
\eqref{spanning-multiplication-operators} give
$$
 P=RL+L^2-LR,\qquad Q=4L^2-3LR-2RL+2R^2.
$$
By \eqref{rsym},
$$
 RQ-QR=R_{(x^2)x}-R_{x(x^2)}.
$$
Using \eqref{spanning-multiplication-operators}, we can write the right side as  
$$
 2PL-3PR-2QL+2QR-2LP+3LQ+2RP-2RQ.
$$
Substituting the expressions for $P$ and $Q$ above and moving all terms
to the left, we obtain
$$
 4RL^2+2RLR-4R^2L-2L^2R-4L^3+4LRL=0,
$$
which is \eqref{spanning-operator-relation}.
\end{proof}
\begin{lemma}\label{metabelian-quartic}
Algebras in $\mathcal V$ satisfy the following identities:
\begin{equation}\label{metabelian}
 [[a,b],[c,d]]=0,
\end{equation}
\begin{equation}\label{powerful identity}
 (ab)[c,d]=a[[c,d],b],
\end{equation}
\begin{equation}\label{powerful identity-2}
    (a[b,c])d-2(ad)[b,c]=0,
\end{equation}
\begin{equation}\label{ident1}
   b[[c,a],d]-b[[d,a],c]+c[[d,a],b]-c[[b,a],d]+d[[b,a],c]-d[[c,a],b]=0.
\end{equation}
\end{lemma}

\begin{proof}
By \eqref{rsym-commutator}, \eqref{half-lie-1},
and \eqref{mainident1},
$$[ab,v]+[a,vb]=[a,v]b+b[a,v]+a(vb)-v(ab)=3[[a,v],b].$$
Setting $D_a=L_a-R_a$, we obtain
\begin{equation}\label{quartic-ad-product}
 D_{ab}=-D_aR_b-3D_bD_a.
\end{equation}
If we set  $s=[b,c]$, then \eqref{rsym-commutator} and the Jacobi identity can be written as $[R_b,R_c]=-R_s$ and $[D_c,D_b]=-D_s$, respectively.
We use \eqref{quartic-ad-product} twice and we obtain
$$D_{as}=D_{(ab)c}-D_{(ac)b}=D_a[R_b,R_c]+9[D_c,D_b]D_a=-D_aR_s-9D_sD_a.$$
On the other hand, from \eqref{quartic-ad-product} follows
$D_{as}=-D_aR_s-3D_sD_a$. Hence $6D_sD_a=0$. Therefore, for arbitrary $d$,
$$D_sD_a(d)=[[b,c],[a,d]]=0.$$

Define
$$
 g(a,b,c,d)=(ab)[c,d]-a[[c,d],b].
$$
Identity \eqref{half-lie-1} yields
$$
 g(a,b,c,d)-g(b,a,c,d)
 =[[a,b],[c,d]]=0.
$$
Note that $g$ is symmetric in its first two variables and alternating
in its last two variables.

Furthermore, we observe that
\begin{equation}\label{quartic-cyclic}
 g(a,b,c,d)+g(a,c,d,b)+g(a,d,b,c)=0.
\end{equation}
To see this,  we consider the cyclic sum over
$b,c,d$. Identities \eqref{rsym} and \eqref{half-lie-1}, together with the Jacobi identity, yield
$$\sum_{\mathrm{cyc}}(ab)[c,d]=\sum_{\mathrm{cyc}}
   \bigl([a,b][c,d]+(ba)[c,d]\bigr)=\sum_{\mathrm{cyc}}
   \bigl([a,b][c,d]+(b[c,d])a+b[a,[c,d]]\bigr)$$
$$=\sum_{\mathrm{cyc}}\bigl([a,b][c,d]+b[a,[c,d]]\bigr)=\sum_{\mathrm{cyc}}\bigl([a,b][c,d]+b[[a,c],d]+b[c,[a,d]]\bigr)$$
$$=\sum_{\mathrm{cyc}}
   \bigl([a,b][c,d]+d[[a,b],c]+c[d,[a,b]]\bigr)=0.$$
Since $\sum_{\mathrm{cyc}}[[b,c],d]=0$, the preceding calculation gives \eqref{quartic-cyclic}.

We next show that $g(x,x,x,y)=0$.
Set $N_a=L_a-2R_a$. For $q=[a,b]$,
\eqref{mainident1} and \eqref{rsym} yield
$$[L_a,L_b]=2L_q-4R_q,\qquad [R_a,R_b]=-R_q,$$ $$[L_a,R_b]+[R_a,L_b]=[L_a,L_b]-L_q.$$
It follows that
\begin{equation}\label{quartic-N-commuting}
 [N_a,N_b]=0.
\end{equation}

If we set $S=2R^2L-2RL^2-RLR+L^2R+2L^3-2LRL$, then $S=0$ by \eqref{spanning-operator-relation}. 
Recall that $P=L_{x^2}$ and by 
\eqref{spanning-multiplication-operators} gives
$P=RL+L^2-LR$. Since $[x,y]=(L-R)(y)$ and $[[x,y],x]=-(L-R)^2y$, \eqref{quartic-N-commuting} gives
$$g(x,x,x,y)=\bigl(P(L-R)+L(L-R)^2\bigr)y=\left(S+\frac13[N_x,N_{x^2}]\right)y=0.$$
We linearize the obtained equality and use symmetry in the first two
variables. Then  we have
\begin{equation}\label{quartic-polarization}
 g(a,b,c,d)+g(a,c,b,d)+g(b,c,a,d)=0.
\end{equation}
For convenience, set
$$
U=g(a,b,c,d),\qquad V=g(a,c,b,d),\qquad W=g(b,c,a,d).
$$
Then $U+V+W=0$. By \eqref{quartic-cyclic} and the symmetries of $g$,
$$
 g(a,d,b,c)=V-U,\qquad
 g(b,d,a,c)=W-U.
$$
We interchange $c$ and $d$ in \eqref{quartic-polarization} and use the last relations yield
$-3U+V+W=0$. Consequently, $4U=0$, and this implies 
\eqref{powerful identity}.
Identities \eqref{rsym-commutator} and \eqref{powerful identity} imply
$$
 (a[b,c])d=a[[b,c],d]+(ad)[b,c]=2(ad)[b,c].
$$
To prove \eqref{ident1}, denote its left-hand side by $E$.
By the Jacobi identity,
$$
E=b[[c,d],a]+c[[d,b],a]+d[[b,c],a].
$$
On the other hand, applying \eqref{powerful identity} and then
\eqref{half-lie-1}, we obtain
$$
E=[b,d][c,a]+[c,b][d,a]+[d,c][b,a]=$$$$
c[[b,d],a]+d[[c,b],a]+b[[d,c],a]-a\bigl([[b,d],c]+[[c,b],d]+[[d,c],b]\bigr)=-E-aJac(b,d,c)=-E.
$$
Hence $2E=0$ and, consequently, $E=0$. 
\end{proof}

We will write the left-normed commutator product as
$$[a_1,\ldots,a_n]=[[\cdots[a_1,a_2]\cdots],a_n]$$
The Jacobi identity and \eqref{metabelian} then give 
\begin{equation}\label{metabelianconseq}
 [a,b,c,d]=[a,b,d,c].
\end{equation}

\begin{lemma}\label{lemma-n-identity-1}
Algebras in $\mathcal{V}$ satisfy
\begin{equation}\label{n-identity-1}
 a_1[a_2,a_3,a_4,\ldots,a_n]
 +a_2[a_3,a_1,a_4,\ldots,a_n]
 +a_3[a_1,a_2,a_4,\ldots,a_n]=0
\end{equation}
for $n\ge3$. Moreover, if $Z=[b_1,\ldots,b_n]$ and $m,n\ge2$, then
\begin{equation}\label{n-identity-2}
 [a_1,\ldots,a_m]Z
 =a_1[Z,a_2,a_3,\ldots,a_m]-a_2[Z,a_1,a_3,\ldots,a_m].
\end{equation}
\end{lemma}

\begin{proof}
For every commutator $y$, \eqref{powerful identity} and
\eqref{powerful identity-2} give
\begin{equation}\label{right-tail}
 (ay)d=2a[y,d].
\end{equation}
The case $n=3$ of \eqref{n-identity-1} is \eqref{half-lie-1}.
Multiply this case of \eqref{n-identity-1}  successively on the right by $a_4,\ldots,a_n$, and using \eqref{right-tail} at each step. Dividing by $2^{n-3}$ gives the first statement for $n>3$.

For $m=2$, identity \eqref{powerful identity} implies \eqref{n-identity-2}. If $U=[a_1,\ldots,a_{m-1}]$ and $m>2$, then
$$
 [U,a_m]Z=U[Z,a_m]-a_m[Z,U]=U[Z,a_m],
$$
since $[Z,U]=0$. Induction on $m$, together with \eqref{metabelianconseq}, completes  the proof.
\end{proof}

Let $F(X)=F_{\mathcal{V}}(X)$ be the free algebra on the countable set
$X=\{x_1,x_2,\ldots\}$. In the following lemmas, let $y$ and $C$ denote left-normed Lie monomials in the generators of length at least two.
For $r\ge0$, set
$$M_r=span\{L_{x_{i_s}}\cdots L_{x_{i_1}}y:0\le s\le r\},
 \qquad M_{-1}=0.$$
Then $M_r\subseteq M_{r+1}$ and $L_xM_r\subseteq M_{r+1}$ for $x\in X$.

\begin{lemma}\label{lemma-n-identity-2}
For $t\ge1$ and $\sigma\in S_t$,
$$
 L_{x_{\sigma(t)}}\cdots L_{x_{\sigma(1)}}y
 \equiv L_{x_t}\cdots L_{x_1}y\pmod{M_{t-1}}.
$$
\end{lemma}

\begin{proof}

The case $t=1$ is trivial, so assume that $t\geq2$.

We first prove a helpful statement. Let $C$ be any left-normed
Lie monomial of length at least two. For generators
$u_1,\ldots,u_r$, define
$$
T_0=y,\qquad T_j=u_jT_{j-1}\quad(1\leq j\leq r).
$$
We show that
\begin{equation}\label{Lie-times-M}
CT_r,\ T_rC\in M_{r+1}
\qquad(r\ge0).
\end{equation}
We prove this by induction on $r$, simultaneously for all choices
of $C$ and $y$.

For $r=0$, identity \eqref{n-identity-2}, applied to $Cy$ and $yC$,
expresses each product as a linear combination of generators
multiplying Lie monomials. Therefore, both $Cy$ and $yC$ belong to $M_1$.

Let $r\ge1$ and assume the statement holds for $r-1$.
By \eqref{powerful identity},
$$
T_rC=(u_rT_{r-1})C
=u_r[C,T_{r-1}]
=u_r(CT_{r-1}-T_{r-1}C).
$$
By the induction hypothesis $CT_{r-1}$ and $T_{r-1}C$ belong to $M_r$. So $T_rC\in M_{r+1}$.

By \eqref{mainident1} we have
$$
CT_r=C(u_rT_{r-1})
=u_r(CT_{r-1})
-4T_{r-1}[C,u_r]+2[C,u_r]T_{r-1}.
$$
Since $[C,u_r]$ is again a left-normed Lie monomial,
the induction hypothesis applies to both $C$ and $[C,u_r]$.
The first term belongs to $M_{r+1}$, and the last two terms
belong to $M_r\subseteq M_{r+1}$.
This proves \eqref{Lie-times-M}.

Now we interchange two adjacent outer multipliers $x_i,x_j\in X$.
If the expression inside them is $T_r$, then
\eqref{mainident1} and \eqref{Lie-times-M}, with $C=[x_i,x_j]$, give
$$
x_i(x_jT_r)-x_j(x_iT_r)
=-4T_r[x_i,x_j]+2[x_i,x_j]T_r\in M_{r+1}.
$$
Since $L_xM_s\subseteq M_{s+1}$ for every $x\in X$,
applying the remaining outer multipliers shows that each adjacent
interchange changes the expression by an element of $M_{t-1}$.
Since adjacent transpositions generate $S_t$, the required congruence holds.
\end{proof}

\begin{lemma}\label{lemma-n-identity-3}
Every left-normed product whose factors are generators except for one factor
$y$ belongs to $M_1$.
\end{lemma}

\begin{proof}
If $V_s=(\cdots(x_1x_2)\cdots)x_s$, repeated use of
\eqref{powerful identity} gives
\begin{equation}\label{left-comb-times-Lie}
 V_sy=x_1[y,x_s,x_{s-1},\ldots,x_2]\in M_1
\end{equation}
for $s\geq1$.
For every generator $z$, we have
$$
yz=zy+[y,z]\in M_1.
$$
Moreover, for generators $x,z_1,\ldots,z_r$, where $r\ge1$,
repeated application of \eqref{right-tail} gives
\begin{equation}\label{successive-right-tails}
(\cdots((xy)z_1)\cdots)z_r
=2^r x[y,z_1,\ldots,z_r]\in M_1.
\end{equation}
\end{proof}

\begin{corollary}\label{corollary-n-identity-4}
For every $\sigma\in S_n$,
$$(\cdots(x_{\sigma(1)}x_{\sigma(2)})\cdots)x_{\sigma(n)}
 \equiv (\cdots(x_1x_2)\cdots)x_n\pmod{M_1}.$$
\end{corollary}

\begin{proof}
For the first two variables, the difference is a left-normed product beginning
with a commutator, so the claim follows from Lemma~\ref{lemma-n-identity-3}.
For any other adjacent pair, \eqref{rsym-commutator} gives
$$(ua)b-(ub)a=u[a,b].$$
After restoring the remaining right factors, the difference again in 
$M_1$ by Lemma~\ref{lemma-n-identity-3}. This proves for the claim all adjacent transpositions.
\end{proof}

\begin{lemma}\label{lemma-n-identity-4}
If $Q_r=x_r(x_{r-1}(\cdots(x_2x_1)\cdots))$ and $r\ge1$, then
$yQ_r,\ Q_ry\in M_r$.
\end{lemma}

\begin{proof}
If $r=1$, use $yx_1=x_1y+[y,x_1]$.
For $r>1$, \eqref{powerful identity} and \eqref{mainident1} give
\begin{align*}
 Q_ry&=x_r(yQ_{r-1}-Q_{r-1}y),\\
 yQ_r&=x_r(yQ_{r-1})-4Q_{r-1}[y,x_r]+2[y,x_r]Q_{r-1}.
\end{align*}
Using induction simultaneously for both products, and for the Lie monomial
$[y,x_r]$ as well as $y$, we see that every displayed term belongs to $M_r$.
\end{proof}

\begin{lemma}\label{lemma-n-identity-5}
For every $\sigma\in S_n$,
$$
 x_{\sigma(n)}(\cdots(x_{\sigma(2)}x_{\sigma(1)})\cdots)
 \equiv x_n(\cdots(x_2x_1)\cdots)\pmod{M_{n-2}}.
$$
\end{lemma}

\begin{proof}
Interchanging $x_{\sigma(2)}$ and $x_{\sigma(1)}$ produces a commutator with
outer left multipliers, and hence an element of $M_{n-2}$.
For an adjacent pair above  $Q_r$, use
$$
 a(bQ_r)-b(aQ_r)=-4Q_r[a,b]+2[a,b]Q_r.
$$
Lemma~\ref{lemma-n-identity-4} implies that the right-hand side belongs to $M_r$.
Multiplying the remaining  outer multipliers puts it in $M_{n-2}\subseteq M_n$. This covers all adjacent transpositions.
\end{proof}

\begin{lemma}\label{lemma-n-identity-6}
For all generators $t,z,x_1,\ldots,x_n$, where $n\ge1$, 
$$
t((\cdots(zx_1)\cdots)x_n)
\equiv z((\cdots(tx_1)\cdots)x_n)\pmod{M_n}.
$$
Moreover, let $V_s$ be a left-normed product of $s\ge2$ generators,
and set $U_q=L_{z_q}\cdots L_{z_1}V_s$, where $q\ge0$ and
$z_1,\ldots,z_q\in X$. Then
\begin{equation}\label{Lie-times-mixed}
 yU_q,\ U_qy\in M_{q+2}.
\end{equation}
\end{lemma}

\begin{proof}
For $n=1$, \eqref{mainident1} yields
$$
t(zx_1)-z(tx_1)
=-2x_1[t,z]+2[[t,z],x_1]\in M_1.
$$

Suppose that $n\ge2$. Set $V_1=x_1$ and
$V_s=V_{s-1}x_s$ for $s\ge2$. First we show  that
$$
V_sy\in M_1,\qquad yV_s\in M_2
$$
for every $s\ge1$ and every left-normed Lie monomial $y$
of length at least two.

The first inclusion follows from
Lemma~\ref{lemma-n-identity-3}.
For the second, we use induction on $s$, simultaneously
for all such $y$. If $s=1$, then
$$
yV_1=x_1y+[y,x_1]\in M_1\subseteq M_2.
$$
For $s\ge2$, using \eqref{degree-3 identity} with
$a=V_{s-1}$, $b=y$, $c=x_s$ gives
$$
yV_s
=-3V_{s-1}[y,x_s]
+2[y,x_s]V_{s-1}
+x_s(V_{s-1}y).
$$
The first term belongs to $M_1$ by
Lemma~\ref{lemma-n-identity-3}.
Since $[y,x_s]$ is again a left-normed Lie monomial,
the second term belongs to $M_2$ by the induction hypothesis.
The last term also belongs to $M_2$, because
$V_{s-1}y\in M_1$ and $x_sM_1\subseteq M_2$.

Using \eqref{degree-3 identity} with
$a=V_n$, $b=z$, $c=t$, we now have
$$
t(V_nz)-z(V_nt)
=3V_n[z,t]-2[z,t]V_n\in M_2.
$$
By Corollary \ref{corollary-n-identity-4} we have
$$V_nz\equiv (\cdots(zx_1)\cdots)x_n\pmod{M_1}, \quad 
V_nt\equiv (\cdots(tx_1)\cdots)x_n\pmod{M_1}.
$$
Multiplying these congruences on the left by $t$ and $z$,
respectively, gives congruences modulo $M_2$. Therefore,
$$
t((\cdots(zx_1)\cdots)x_n)
\equiv z((\cdots(tx_1)\cdots)x_n)\pmod{M_2}.
$$
Since $M_2\subseteq M_n$ for $n\ge2$, the statement follows.

For \eqref{Lie-times-mixed}, the case $q=0$ is exactly the pair
of inclusions proved above. For $q\ge1$, we use  identities
\eqref{powerful identity} and \eqref{mainident1}:
\begin{align*}
 U_qy&=z_q(yU_{q-1}-U_{q-1}y),\\
 yU_q&=z_q(yU_{q-1})
       -4U_{q-1}[y,z_q]+2[y,z_q]U_{q-1}.
\end{align*}
Induction on $q$, applied simultaneously to all generator commutator monomials $y$,
puts the expressions in parentheses and the last two products in
$M_{q+1}$. Left multiplication by $z_q$ puts the remaining terms
in $M_{q+2}$. This proves \eqref{Lie-times-mixed}.
\end{proof}
\begin{corollary}\label{corollary-n-identity-7}
Let $n\ge2$ and $0\le k\le n-2$. For every $\sigma\in S_n$,
$$L_{x_{\sigma(1)}}\cdots L_{x_{\sigma(k)}}
\bigl((\cdots(x_{\sigma(k+1)}x_{\sigma(k+2)})\cdots)x_{\sigma(n)}\bigr)\equiv
 L_{x_1}\cdots L_{x_k}\bigl((\cdots(x_{k+1}x_{k+2})\cdots)x_n\bigr)\pmod{M_{n-2}}.$$
\end{corollary}

\begin{proof}
Define $U_q=L_{z_q}\cdots L_{z_1}V_s$, where $s\ge2$ as in Lemma \eqref{lemma-n-identity-6}. For generators $x_i,x_j$, by \eqref{mainident1}  and \eqref{Lie-times-mixed} we have
$$
 x_i(x_jU_q)-x_j(x_iU_q)
 =-4U_q[x_i,x_j]+2[x_i,x_j]U_q\in M_{q+2}.
$$
It proves that two adjacent outer multipliers can also be interchanged and hence any further outer multipliers increase the
filtration bound by at most their number. The other remaining cases are immediate consequences of Corollary~\ref{corollary-n-identity-4}
and Lemmas~\ref{lemma-n-identity-5}, \ref{lemma-n-identity-6}.
\end{proof}

\section{A basis and multilinear dimension}
We order the generators by their indices as $x_1<x_2<\cdots$. The elements of $Y_1$ and $Y_2$ have the respective forms
\begin{equation}\label{spanning-unrestricted}
\begin{aligned}
 &L_{x_{i_t}}\cdots L_{x_{i_1}}
       [x_j,x_{k_1},\ldots,x_{k_s}], &&t\ge0,\quad s\ge1,\\
 &L_{x_{j_t}}\cdots L_{x_{j_1}}
       (((x_{i_1}x_{i_2})\cdots)x_{i_s}), &&t\ge0,\quad s\ge2.
\end{aligned}
\end{equation}
The indices in $Y_1$  satisfy
$$
 j>k_1\le k_2\le\cdots\le k_s,
 \qquad k_1\le i_1\le\cdots\le i_t,
$$
and there is no triple $i_p>j>k_q$ with $1\le p\le t$ and
$2\le q\le s$. The indices in $Y_2$ satisfy
$$
 j_t\ge\cdots\ge j_1\ge i_1\ge\cdots\ge i_s.
$$
We make the restriction $s\ge2$ to avoid two descriptions of the same right-normed monomial. We separately include generators of degree one.

Write $\widetilde Y_1,\widetilde Y_2$ for the same families with all
ordering without restrictions above. Thus
$$
 M=\bigcup_{r\ge0}M_r=span\widetilde Y_1.
$$
Let $P_n(\mathcal V)$ denote  the multilinear component in
$x_1,\ldots,x_n$ of the free algebra
$F_{\mathcal V}(X)$.
\begin{theorem}\label{basis-theorem}
For every $n\ge2$, the multilinear elements of $Y_1\cup Y_2$
in $x_1,\ldots,x_n$ form a basis of $P_n(\mathcal V)$.
\end{theorem}
\begin{proof}
Taking the multilinear component in Lemma~\ref{spanning-Y},
we get that the multilinear elements of $Y_1\cup Y_2$
span $P_n(\mathcal V)$. Lemma~\ref{independence-Y} provides the independence of them.
\end{proof}
\begin{lemma}\label{spanning-ordering}
The space $M$ is spanned by $Y_1$. Moreover,
$$
 span(X\cup\widetilde Y_1\cup\widetilde Y_2)
 =span(X\cup Y_1\cup Y_2).
$$
\end{lemma}

\begin{proof}
We prove $M_t\subseteq span Y_1$ by induction on $t$. The
Jacobi identity and \eqref{metabelianconseq} allow us to order the commutator part with its least variable in the second position. This proves the case $t=0$.
For $t>0$, Lemma~\ref{lemma-n-identity-2} allows us order the outer
multipliers modulo $M_{t-1}$. If the least variable $m$ occurs
only among these multipliers, use \eqref{n-identity-1} in the form
$$
 m[b,c,E]=c[b,m,E]-b[c,m,E]
$$
to bring it into the second position of the commutator part.

To remove a forbidden triple $a>b>c$, apply \eqref{ident1},
extended along the tail $E$ by \eqref{right-tail}:
$$
\begin{aligned}
 a[b,m,c,E]={}&a[c,m,b,E]-b[c,m,a,E]\\
 &+b[a,m,c,E]-c[a,m,b,E]+c[b,m,a,E].
\end{aligned}
$$
After restoring and ordering the remaining outer multipliers,
each nonzero term either belongs to $M_{t-1}$ or has a strictly
smaller pair, consisting of the increasing tuple of outer indices 
and the first commutator index, in lexicographic order.
At a fixed multidegree, there are only finitely many such pairs,
so the reduction terminates in $Y_1$.

Finally, Corollary~\ref{corollary-n-identity-7} allows us to order every
element of $\widetilde Y_2$ into an element of $Y_2$ modulo
$M=span\,Y_1$, and proving the second statement.
\end{proof}

\begin{lemma}\label{spanning-Y}
The set $X\cup Y_1\cup Y_2$ spans $F(X)$.
\end{lemma}

\begin{proof}
Put
$$
 \mathcal{W}=span(X\cup\widetilde Y_1\cup\widetilde Y_2),
 \qquad \mathcal{W}_n=\mathcal{W}\cap F(X)_n.
$$
By Lemma~\ref{spanning-ordering},
$\mathcal{W}=span(X\cup Y_1\cup Y_2)$.
Left multiplication by a generator preserves $\mathcal{W}$: it adds an
outer left multiplier to an unrestricted normal form, and the
product of two generators belongs to $\widetilde Y_2$.
Thus the identities \eqref{spanning-multiplication-operators} reduces the proof to showing that
$$ux\in\mathcal{W}
 \qquad(x\in X,\quad u\in\widetilde Y_1\cup\widetilde Y_2).$$
We use induction on the total degree $n$.
Every monomial of degree at most three already belongs to $\mathcal{W}$.
Fix $n\ge4$ and assume that every monomial of smaller degree
belongs to $\mathcal{W}$. Hence, by linearity, it suffices to treat the
displayed products with $\deg u=n-1$.
Throughout the rest of the proof, $f\equiv g$ means
$f-g\in\mathcal{W}_n$, and every congruence has total degree $n$.
In particular, the induction hypothesis and $x\mathcal{W}\subseteq\mathcal{W}$ yield 
\begin{equation}\label{spanning-left-generator}
 xa\equiv0
 \qquad(x\in X,\quad a\in F(X)_{n-1}).
\end{equation}

We first establish the reductions needed for the products $ux$.
For homogeneous $a,b,c$ whose degrees sum to $n$,  the identity \eqref{powerful identity} and  \eqref{spanning-left-generator} give 
\begin{equation}\label{spanning-right-commutator}
 a[b,c]\equiv0.
\end{equation}
Consequently, right-symmetry gives
\begin{equation}\label{spanning-exchanges}
 (ab)c\equiv(ac)b,\qquad a(bc)\equiv a(cb).
\end{equation}
Applying \eqref{degree-3 identity} to $b, a, x$, where $x\in X$,
and using \eqref{spanning-left-generator} and
\eqref{spanning-right-commutator}, we also obtain
\begin{equation}\label{spanning-transfer}
 a(xb)\equiv a(bx)\equiv2[a,x]b.
\end{equation}

Suppose that $T$ is obtained from a commutator by successive
left or right multiplication by generators. Then
\begin{equation}\label{spanning-context}
 aT\equiv0.
\end{equation}
We prove the statment for every left factor $a$ and
every initial commutator by induction on the number
of these multiplications.
The case of no multiplications is
\eqref{spanning-right-commutator}.
If $T=xS$, then $aT\equiv2[a,x]S$, and the induction hypothesis
applies to $S$. If $T=Sx$, first replace $a(Sx)$ by $a(xS)$
using \eqref{spanning-exchanges}.

Let $U$ be obtained from a generator by successive left or
right multiplications by generators. Replace every left
multiplication $xS$ by $Sx$, keeping the order of the steps,
and denote the resulting left-normed product by $V$. Then
\begin{equation}\label{spanning-comb-orientation}
 aU\equiv aV.
\end{equation}
To see this, let $P$ denote the composition of the remaining
left and right multiplications by generators
(possibly the identity operator). Each replacement gives
$$
 aP(xS)-aP(Sx)=aP([x,S])\equiv0
$$
by \eqref{spanning-context}.
Repeating this for every left multiplication yields
$aU\equiv aV$.

We also have
\begin{equation}\label{spanning-M-products}
 af\equiv fa\equiv0
 \qquad(f\in M).
\end{equation}
The first congruence follows from \eqref{spanning-context}.
For the second, first take $f=y$, a left-normed Lie
monomial in the generators. By the degree induction hypothesis, $a$ is a linear
combination of unrestricted normal forms and generators.
For a generator $z$, we have $yz=zy+[y,z]\in M$.
For the other forms, \eqref{Lie-times-M} and
\eqref{Lie-times-mixed} give $ya\in M\subseteq\mathcal{W}$.
If $f=zU\in\widetilde Y_1$, where $z\in X$ and $U\in M$,
then
$$
 fa=(zU)a\equiv(za)U\equiv0
$$
by \eqref{spanning-exchanges} and \eqref{spanning-context}.
This proves \eqref{spanning-M-products}.

We now consider $ux$, with $x\in X$ and
$u\in\widetilde Y_1\cup\widetilde Y_2$.
If $u\in\widetilde Y_1$, then $ux\equiv0$ by
\eqref{spanning-M-products}.
If $u\in\widetilde Y_2$ has no outer left multiplier, then
$ux\in\widetilde Y_2$.
Otherwise, write $u=aU$, where $a\in X$, and let $V$ be the
left-normed generator product associated with $U$ in
\eqref{spanning-comb-orientation}. Then
$$
 ux=(aU)x\equiv(ax)U\equiv(ax)V\equiv(aV)x.
$$
Thus it remains to prove that
\begin{equation}\label{spanning-T}
 T_n(x_1,\ldots,x_n)
 =\bigl(x_1(((x_2x_3)\cdots)x_{n-1})\bigr)x_n
 \in\mathcal{W}_n.
\end{equation}

The class of $T_n$ modulo $\mathcal W_n$ is symmetric
in all its variables. By
Corollary~\ref{corollary-n-identity-7} and
\eqref{spanning-M-products}, it is invariant under
permutations of the first $n-1$ variables.
Writing $V=(\cdots(x_2x_3)\cdots)x_{n-1}$, we also have
$$
 (x_1V)x_n\equiv(x_1x_n)V
 \equiv(x_nx_1)V\equiv(x_nV)x_1.
$$
Here the middle congruence follows from
$[x_1,x_n]V\equiv0$.
These permutations generate $S_n$.

Fix a generator $x$, put $L=L_x$, $R=R_x$, and write
$T=T_n(x,\ldots,x)$.
For $0\le p\le n-2$, we claim that
\begin{equation}\label{spanning-recurrence}
 R^pLf\equiv c_pT
 \qquad
 \bigl(f=S_1\cdots S_{n-p-2}x,\quad S_i\in\{L,R\}\bigr),
\end{equation}
where an empty operator word is allowed and
$$
 c_0=0,\qquad c_1=1,\qquad
 c_p=\frac{3c_{p-1}-c_{p-2}}2\quad(p\ge2).
$$
The case $p=0$ is \eqref{spanning-left-generator}.
For $p=1$, set $V=R^{n-3}x$. Then we have
\begin{equation}\label{spanning-RL}
 RLf=(xf)x\equiv(xx)f\equiv(xx)V\equiv(xV)x=T,
\end{equation}
by \eqref{spanning-exchanges} and
\eqref{spanning-comb-orientation}.

For $p\ge2$, multiply the exact operator identity
\eqref{spanning-operator-relation} by $R^{p-2}$ on the left
and apply it to $f$. Then we get
$$2R^pLf=2R^{p-1}L(Lf)+R^{p-1}L(Rf)-R^{p-2}L(LRf)-2R^{p-2}L(L^2f)+2R^{p-2}L(RLf).$$
Each argument on the right is again an operator word applied to
$x$, of the required degree. So induction on $p$ therefore implies
$$
 2R^pLf\equiv(3c_{p-1}-c_{p-2})T,
$$
and this proves \eqref{spanning-recurrence}.

One can have the solution to the recurrence relation
$$
 c_p=\frac{2^p-1}{2^{p-1}}\qquad(p\ge1).
$$
Set $p=n-2$ and $f=x$. Since $Lx=Rx$,
$$
 R^{n-2}Lx=R^{n-1}x\in\widetilde Y_2.
$$
Hence
$$
 \frac{2^{n-2}-1}{2^{n-3}}T\equiv0.
$$
The coefficient is nonzero in the field of characteristic zero, so
$T_n(x,\ldots,x)\in\mathcal{W}_n$.

Completely linearizing it, we have
$$
 \sum_{\pi\in S_n}
 T_n(x_{\pi(1)},\ldots,x_{\pi(n)})\in\mathcal{W}_n,
$$
because the complete linearization of an unrestricted normal form
is a sum of unrestricted normal forms.
By the symmetry established above, this implies
$n!T_n\in\mathcal{W}_n$, and therefore $T_n\in\mathcal{W}_n$.

Thus every required product $ux$ of degree $n$ belongs to $\mathcal{W}$.
Together with $x\mathcal{W}\subseteq\mathcal{W}$ and the identities 
\eqref{spanning-multiplication-operators}, this completes the induction.
Lemma~\ref{spanning-ordering} gives the stated spanning set.
\end{proof}

\begin{lemma}\label{independence-Y}
For every $n\ge2$, the multilinear elements of $Y_1\cup Y_2$
in $x_1,\ldots,x_n$ are linearly independent modulo
$\operatorname{Id}(B)$.
\end{lemma}

\begin{proof}
Write $w_{j,I}$ for the multilinear $Y_1$-element whose first
commutator index is $j$ and whose outer index set is $I$.
Its second commutator index is $1$, and the remaining indices form
the increasing commutator tail. Let $Q_r$ denote the ordered $Y_2$-element with $n-r$ outer left multipliers, where $2\le r\le n$.

Suppose that
$$
 f=\sum_{(j,I)}a_{j,I}w_{j,I}
       +\sum_{r=2}^{n}b_rQ_r
 \in\operatorname{Id}(B).
$$
We prove that all coefficients vanish using only substitutions
of the basis elements $e_1,e_2$.

For $S\subseteq\{1,\ldots,n\}$, let $\varphi_S$ be the substitution
$$
 x_i=
 \begin{cases}
 e_2,&i\in S,\\
 e_1,&i\notin S.
 \end{cases}
$$
Put
$$
 E_k=L_{e_2}^{\,k-1}(e_2)
 =
 \begin{cases}
 e_1,&k\text{ even},\\
 e_2,&k\text{ odd}.
 \end{cases}
$$

A nonzero commutator evaluation requires exactly one of its
first two variables to be $e_2$ and every tail variable to be
$e_1$. Consequently,
$$
 \varphi_S(w_{j,I})=0\qquad\text{if }|S|>|I|+1.
$$
If $|S|=|I|+1=k$, then
\begin{equation}\label{independence-Y1-evaluation}
 \varphi_S(w_{j,I})=
 \begin{cases}
 E_k,&S=I\cup\{j\},\\
 -E_k,&S=I\cup\{1\},\\
 0,&\text{otherwise}.
 \end{cases}
\end{equation}

Similarly, a nonzero left-normed product contains at most
two occurrences of $e_2$. Thus
$$\varphi_S(Q_r)=0\qquad\text{if }|S|>n-r+2.$$
When $|S|=n-r+2=k$, all outer variables must be $e_2$,
and the two variables taking the value $e_2$ in the
left-normed product must be $x_r$ and some $x_i$ with $i<r$. Hence
\begin{equation}\label{independence-Y2-evaluation}
 \varphi_S(Q_r)=
 \begin{cases}
 2^{i-1}E_k,
   &S=\{r,r+1,\ldots,n\}\cup\{i\},\quad 1\le i<r,\\
 0,&\text{otherwise}.
 \end{cases}
\end{equation}
Indeed, after the product reaches the second $e_2$,
its value is $e_1$, and each of the remaining $i-1$
right factors $e_1$ contributes a factor of $2$.

We use descending induction on $k=n,n-1,\ldots,1$.
At the beginning of the $k$th step, assume that
$$
 a_{j,I}=0 \quad\text{if } |I|+1>k,
 \qquad
 b_r=0 \quad\text{if } n-r+2>k.
$$
For $k=n$, this assumption holds automatically.
Our aim is to prove that all coefficients with $|I|=k-1$
vanish and, if $k\ge2$, that $b_{n-k+2}=0$.

Consider substitutions $\varphi_S$ with $|S|=k$.
Terms with $|I|+1>k$ or $n-r+2>k$ have zero coefficients
by the induction hypothesis. By the preceding observations, terms with $|I|+1<k$ or
$n-r+2<k$ vanish under these substitutions.
Thus only the $Y_1$-terms with $|I|=k-1$ and, if $k\ge2$,
the $Y_2$-term with $r=n-k+2$ can remain nonzero under
these substitutions.

For $k\ge2$, sum the equations $\varphi_S(f)=0$ over all
$k$-element subsets $S$. By \eqref{independence-Y1-evaluation},
each remaining term $a_{j,I}w_{j,I}$ becomes $a_{j,I}E_k$
when $S=I\cup\{j\}$ and $-a_{j,I}E_k$ when
$S=I\cup\{1\}$, and is zero for every other
$k$-element subset $S$.
These two subsets are distinct and both occur in the sum,
so the two values cancel. Hence the sum of all substituted
$Y_1$-terms is zero. Therefore, for $r=n-k+2$,
$$
 (1+2+\cdots+2^{r-2})b_rE_k=0.
$$
Since the field has characteristic zero and $E_k\ne0$,
we obtain $b_r=0$. In particular, the initial step $k=n$
gives $b_2=0$ and there are no $Y_1$-terms with $|I|=n-1$.
As $k$ decreases from $n$ to $2$, we successively obtain
$b_2=b_3=\cdots=b_n=0$. Thus, at the final step $k=1$,
no $Y_2$-term remains, and only the $Y_1$-terms with
no outer multipliers need to be considered.

Now we fix $k$. It remains to eliminate the coefficients
$a_{j,I}$ with $|I|=k-1$. For this purpose, use a second
descending induction on $m=n,n-1,\ldots,2$.
At the beginning of the $m$th step, assume that
$$
 a_{j,I}=0
 \quad\text{if } |I|=k-1
 \text{ and } \max(I\cup\{j\})>m.
$$
Clearly, the second induction hypothesis holds for $m=n$.
We prove that all coefficients with $|I|=k-1$ and
$\max(I\cup\{j\})=m$ also vanish, considering first
$m\in I$ and then $j=m$.

Suppose first that $m\in I$, so $j<m$. If some
$h\in\{2,\ldots,j-1\}$ is not in $I$, then $h$
 belongs to the commutator tail, giving the forbidden triple
$m>j>h$. Hence
$$
 \{2,\ldots,j-1\}\subseteq I,
 \qquad
 j=\min\bigl(\{2,\ldots,m\}\setminus I\bigr).
$$
Thus, for this outer set $I$, at most one initial index
can still have a nonzero coefficient.

Apply the substitution $\varphi_{I\cup\{1\}}$.
Since this set contains $1$, formula
\eqref{independence-Y1-evaluation} shows that only terms
with outer set $I$ can remain nonzero under this substitution.
Each such term $a_{j,I}w_{j,I}$ becomes $-a_{j,I}E_k$.
By the preceding uniqueness and the induction hypothesis,
the equation reduces to
$$
 -a_{j,I}E_k=0.
$$
Therefore all coefficients under consideration with
$m\in I$ vanish.

Now consider a remaining coefficient $a_{m,I}$.
Apply $\varphi_{I\cup\{m\}}$. This set does not contain $1$,
so a contributing term $w_{j',I'}$ must satisfy
$$
 I'\cup\{j'\}=I\cup\{m\}.
$$
If $j'=m$, then $I'=I$, giving the chosen term.
If $j'\ne m$, then $m\in I'$, and its coefficient has
just been proved zero. Consequently,
$$
 a_{m,I}E_k=0,
$$
so every remaining coefficient with first commutator index $m$
also vanishes.

We have therefore proved that
$$
 a_{j,I}=0
 \quad\text{if } |I|=k-1
 \text{ and } \max(I\cup\{j\})=m.
$$
We can now repeat the argument with $m-1$ in place of $m$.
After reaching $m=2$, every coefficient with $|I|=k-1$
is zero.

The $k$th step is therefore complete, and the hypothesis
for the next step $k-1$ follows. We continue down to $k=1$ and conclude that all coefficients in $f$ vanish. This completes the proof.
\end{proof}
\begin{corollary}\label{finite-basis-identities}
The $T$-ideal $T(B)$ is generated by the right-symmetric
identity \eqref{rsym} and the identities
\eqref{half-lie-1}, \eqref{irred-(2,1)}.
Consequently, $\mathcal V=\operatorname{var}(B)$.
\end{corollary}

\begin{proof}
Let $J$ be the $T$-ideal generated by these three identities.
By Proposition~\ref{B-identities} and the right-symmetry
of $B$, we have $J\subseteq T(B)$.

Let $f$ be a multilinear identity of $B$ of degree $n\ge2$.
By Lemma~\ref{spanning-Y},
$$
 f\equiv\sum_w\alpha_w w\pmod J,
$$
where $w$ ranges over the multilinear elements of
$Y_1\cup Y_2$ in $x_1,\ldots,x_n$.
Since $J\subseteq T(B)$, the polynomial
$\sum_w\alpha_w w$ is also an identity of $B$.
Lemma~\ref{independence-Y} gives $\alpha_w=0$ for every $w$,
and hence $f\in J$.

The algebra $B$ has no nonzero identities of degree one.
Since the  field has characteristic zero, every
polynomial identity follows from multilinear identities.
Therefore $T(B)=J$, and
$\mathcal V=\operatorname{var}(B)$.
\end{proof}
\begin{corollary}\label{multilinear-dimension}
For every $n\ge2$,
$$
c_n(\mathcal V)=\dim_{\mathbb F}P_n(\mathcal V)=2^n-2.
$$
\end{corollary}
\begin{proof}
By Theorem~\ref{basis-theorem}, it suffices to count the
multilinear elements of $Y_1\cup Y_2$.

Let $Y_{1,t}(n)$ denote the multilinear $Y_1$-elements
with exactly $t$ outer multipliers, where $0\le t\le n-2$.
Each such element is uniquely determined by its first
commutator index $j$ and its outer index set
$$
I\subseteq\{2,\ldots,n\}\setminus\{j\},
\qquad |I|=t.
$$
Its second commutator index is $1$, and the remaining
$n-t-2$ indices form the increasing commutator tail.

First suppose that $j=\max(I\cup\{j\})$.
There is no forbidden triple, since all outer indices
are smaller than $j$. Choosing the $(t+1)$-element set
$I\cup\{j\}$ determines $j$ as its largest element
and $I$ as the set of its remaining elements.
This gives
$$
\binom{n-1}{t+1}
$$
elements.

If $I$ contains an index greater than $j$, the
forbidden-triple condition forces all indices
$2,\ldots,j-1$ to belong to $I$. Thus $j$ is smaller
than every tail index.

To count the elements with this property, choose
$n-t-1$ indices from $\{2,\ldots,n\}$.
Take the smallest chosen index as $j$ and arrange
the other chosen indices increasingly in the
commutator tail. The $t$ unchosen indices form $I$.
Every such choice is admissible, so this gives
$$
\binom{n-1}{n-t-1}=\binom{n-1}{t}
$$
elements.

The two counts include all admissible elements,
but they overlap precisely when
$$
I=\{2,\ldots,j-1\},
$$
that is, when $j=t+2$ and $I=\{2,\ldots,t+1\}$.
Hence exactly one element is counted twice, and
$$
|Y_{1,t}(n)|
=\binom{n-1}{t+1}+\binom{n-1}{t}-1
=\binom n{t+1}-1.
$$
Summing over $t$ yields
$$
|Y_1(n)|
=\sum_{t=0}^{n-2}\left(\binom n{t+1}-1\right)
=2^n-n-1.
$$

An ordered multilinear $Y_2$-element is uniquely determined
by the length of its inner left-normed product, which
ranges from $2$ to $n$. Thus $|Y_2(n)|=n-1$, and
$$
c_n(\mathcal V)
=|Y_1(n)|+|Y_2(n)|
=(2^n-n-1)+(n-1)
=2^n-2.
$$
\end{proof}

\section{$S_n$-module structures}\label{sn-gen:section}
In this section, we study the $S_n$-module structure of $P_n(\mathcal V)$ and give its complete decomposition into irreducible $S_n$-modules. The symmetric group $S_n$ acts on multilinear polynomials by
$$
 \sigma f(x_1,\ldots,x_n)
 =
 f(x_{\sigma(1)},\ldots,x_{\sigma(n)}).
$$

Denote by $S^\lambda$ the Specht $\mathbb F S_n$-module which is irreducible $\mathbb F S_n$-module corresponding to a partition
$\lambda\vdash n$. We use the standard Young-symmetrizer
theorem and the hook-length formula, for example see \cite{ FultonHarris} or \cite{Sagan}.

Fix $n\ge2$. For $-1\le k\le n-2$, set
$$
M_{n,k}=M_k\cap P_n(\mathcal V),\qquad
\overline M_{n,k}=M_{n,k}/M_{n,k-1}\quad(k\ge0),
\qquad
\mathcal M_n=M_{n,n-2}.
$$
Then $M_{n,-1}=0$, and $M_{n,k}$ is spanned by the multilinear
words
$$
L_{x_{i_t}}\cdots L_{x_{i_1}}
[x_{j_1},\ldots,x_{j_{n-t}}],
\qquad 0\le t\le k.
$$
Since the commutator part has length at least two, the filtration is
$$
0=M_{n,-1}\subseteq M_{n,0}\subseteq\cdots
\subseteq M_{n,n-2}=\mathcal M_n.
$$

\begin{lemma}\label{sn-gen:filtration}
$M_{n,k}$ is an $S_n$-submodule. The classes of the multilinear
$Y_1$-elements with  $k$ outer multipliers is a basis of
$\overline M_{n,k}$, and
\begin{equation}\label{sn-gen:layer-dimension}
\dim_{\mathbb F}\overline M_{n,k}=\binom n{k+1}-1.
\end{equation}
Moreover, we have the following $S_n$-module isomorphism
\begin{equation}\label{sn-gen:trivial-quotient}
P_n(\mathcal V)/\mathcal M_n
\cong (n-1)S^{(n)}.
\end{equation}
\end{lemma}

\begin{proof}

Permuting the variables changes only their indices, so a word
with at most $k$ outer multipliers remains in $M_{n,k}$.
Thus $M_{n,k}$ is an $S_n$-submodule.

The proof of Lemma~\ref{spanning-ordering} expresses each such
word as a linear combination of ordered $Y_1$-words with at
most $k$ outer multipliers. Indeed, the identities used there either
preserve the number of outer multipliers or produce terms
with fewer outer multipliers.

Let $Y_{1,t}(n)$ denote the multilinear $Y_1$-words with
exactly $t$ outer multipliers. Then
$$
 M_{n,k}
 =
 \operatorname{span}_{\mathbb F}
 \bigcup_{t=0}^{k}Y_{1,t}(n).
$$
By Theorem~\ref{basis-theorem}, these words are linearly
independent and therefore give a basis of $M_{n,k}$.
Those with fewer than $k$ outer multipliers form a basis
of $M_{n,k-1}$. Hence the classes of the words with exactly
$k$ outer multipliers form a basis of $\overline M_{n,k}$.

The proof of
Corollary~\ref{multilinear-dimension} gives
$$
|Y_{1,k}(n)|
=\binom{n-1}{k+1}+\binom{n-1}{k}-1
=\binom n{k+1}-1,
$$
which proves \eqref{sn-gen:layer-dimension}.

For $2\le s\le n$, recall
$$
Q_s=L_{x_n}\cdots L_{x_{s+1}}
\bigl((\cdots(x_sx_{s-1})\cdots)x_1\bigr).
$$
By Theorem~\ref{basis-theorem}, the classes
$Q_2+\mathcal M_n,\ldots,Q_n+\mathcal M_n$ form a basis
of $P_n(\mathcal V)/\mathcal M_n$.
In multilinear degree $n$, the  space $M_n$ satisfies
$M_n\cap P_n(\mathcal V)=\mathcal M_n$, because every commutator part
has length at least two. Thus
Corollary~\ref{corollary-n-identity-7} yields
$$
\sigma Q_s-Q_s\in\mathcal M_n
\qquad(\sigma\in S_n).
$$
Every basis class is therefore fixed by $S_n$, and we get \eqref{sn-gen:trivial-quotient}.
\end{proof}

We now construct the irreducible submodules of each quotient
$\overline M_{n,k}$. Let $n\ge2$ be the total degree, and fix an integer
$k$ satisfying $0\le k\le n-2$.
Choose the ordered multilinear $Y_1$-word
\begin{equation}\label{nv:seed}
 u_{n,k}
 =
 L_{x_{n-1}}\cdots L_{x_{n-k}}
 [x_n,x_1,\ldots,x_{n-k-1}].
\end{equation}
The commutator part has length $n-k$, and the remaining
$k$ variables are outer multipliers. We note that by its construction $u_{n,k}$ is one of the prescribed ordered
basis words.

Choose an integer $j$ satisfying
\begin{equation}\label{nv:range}
 1\le j\le\min(k+1,n-k-1),
\end{equation}
and let $T=T_{n,j}$ be the Young tableau of shape
$(n-j,j)$ filled consecutively along its two rows:
\begin{equation}\label{nv:tableau}
 T=T_{n,j}
 =
 \begin{array}{cccccc}
 1&2&\cdots&j&\cdots&n-j\\
 n-j+1&n-j+2&\cdots&n&&
 \end{array}.
\end{equation}
For $1\le i\le j$, the two-entry column numbered $i$
contains $i$ above $n-j+i$.

We now describe the positions of the variables of
$u_{n,k}$ relative to the two rows of $T$.
Set
\begin{equation}\label{nv:parameters}
 r=n-j,\qquad \ell=n-k-1,\qquad p=k-j+1.
\end{equation}
Thus
$$
 r=\ell+p,\qquad k=p+j-1,\qquad
 j\le\ell\le r,\qquad p\ge0.
$$

The following Young tableau indicates the sizes
of the three blocks. 

\begin{center}
\begin{tikzpicture}[
  every node/.style={font=\normalsize},
  brace/.style={
    decorate,
    decoration={brace,amplitude=5pt}
  }
]

\matrix (T) [
  matrix of math nodes,
  nodes={
    minimum height=4.5mm,
    inner xsep=2pt,
    inner ysep=0pt
  },
  column sep=2mm,
  row sep=0mm
] {
  1 & \cdots & j-1 & j & j+1 & \cdots
    & \ell & \ell+1 & \cdots & r \\
  r+1 & \cdots & n-1 & n
    & & & & & & \\
};

\draw[brace]
  ([yshift=3pt]T-1-1.north west) --
  ([yshift=3pt]T-1-7.north east)
  node[midway,above=8pt] {$\ell$};

\draw[brace]
  ([yshift=3pt]T-1-8.north west) --
  ([yshift=3pt]T-1-10.north east)
  node[midway,above=8pt] {$p$};

\draw[
  decorate,
  decoration={brace,amplitude=5pt,mirror}
]
  ([yshift=-3pt]T-2-1.south west) --
  ([yshift=-3pt]T-2-3.south east)
  node[midway,below=8pt] {$j-1$};

\end{tikzpicture}
\end{center}

The first-row variables $x_1,\ldots,x_\ell$ belong
to the commutator part, whereas the remaining $p$
first-row variables $x_{\ell+1},\ldots,x_r$
are outer multipliers.
Among the second-row variables, $x_n$ is the initial
entry of the commutator part, and the remaining $j-1$
variables $x_{r+1},\ldots,x_{n-1}$ are outer multipliers.

Using \eqref{nv:parameters}, we  rewrite
$u_{n,k}$ as
\begin{equation}\label{nv:seed-split}
 u_{n,k}
 =
 \underbrace{L_{x_{n-1}}\cdots L_{x_{r+1}}}
 _{j-1\ \text{second-row variables}}
 \underbrace{L_{x_r}\cdots L_{x_{\ell+1}}}
 _{p\ \text{first-row variables}}
 [x_n,x_1,\ldots,x_\ell].
\end{equation}
Blocks containing zero multipliers are omitted.

Let
$$
 \mathcal R(T)
 =
 S_{\{1,\ldots,r\}}\times S_{\{r+1,\ldots,n\}},
 \qquad
 s_i=(i,r+i)\quad(1\le i\le j).
$$
Define
$$
 R_T=\sum_{\rho\in\mathcal R(T)}\rho,\qquad
 C_T=\prod_{i=1}^{j}(1-s_i),\qquad
 e_T=R_TC_T.
$$
Set
\begin{equation}\label{nv:generator}
 f_{n,k,j}=e_Tu_{n,k}=R_TC_Tu_{n,k},
\end{equation}
where $C_T$ acts on $u_{n,k}$ first, then $R_T$. 
Since $u_{n,k}\in M_{n,k}$ and $M_{n,k}$ is invariant under
the action of $S_n$, we have $f_{n,k,j}\in M_{n,k}$.

\begin{lemma}\label{sn-gen:nonvanishing}
Under \eqref{nv:range}, make the substitution in $B$ given by
$$
E:\quad
x_1=\cdots=x_\ell=e_1,\qquad
x_{\ell+1}=\cdots=x_n=e_2.
$$
Then
\begin{equation}\label{sn-gen:evaluation}
 f_{n,k,j}(E)\ne0.
\end{equation}
Moreover,
$$
 f_{n,k,j}+M_{n,k-1}\ne0
 \qquad\text{in }\overline M_{n,k},
$$
where $M_{n,-1}=0$.
\end{lemma}

\begin{proof}
We first check the column alternation. Suppose $j\ge2$. Fix a row permutation
$\rho\in\mathcal R(T)$ and a nonidentity column permutation
$\tau$ occurring in the expansion of $C_T$.
If $\tau$ contains a transposition $s_i$ with $i\ge2$,
then the original occurrence of $x_i$ in the commutator part
is replaced by
$x_{\rho(r+i)}$ in $\rho\tau u_{n,k}$ .
Since $2\le i\le j\le\ell$, this is a tail position.
Moreover, $\rho(r+i)$ belongs to the second row, so its
value under $E$ is $e_2$.
The preceding commutator belongs to $[B,B]=\mathbb{F} \langle e_2\rangle$, and consequently, this tail entry makes the commutator part zero.

The only remaining possibility is $\tau=s_1$.
Since $j\ge2$, we have $r+1\ne n$.
Consequently, the initial bracket in
$\rho s_1u_{n,k}$ is
$$
 [x_{\rho(n)},x_{\rho(r+1)}].
$$
Both variables are in the second row and receive
$e_2$, so this bracket is zero.
Hence every nonidentity column term vanishes after
every row permutation, and consequently, for $j\ge2$, we have
$$f_{n,k,j}(E)=\sum_{\rho\in\mathcal R(T)}(\rho u_{n,k})(E).$$
Now assume $j=1$. Then the only column transposition is
$s_1=(1,n)$. It exchanges the first two entries of the
commutator part and leaves the outer multipliers unchanged.
Therefore
$$
 s_1u_{n,k}=-u_{n,k},
 \qquad
 C_Tu_{n,k}=2u_{n,k},
$$
and hence
$$
 f_{n,k,1}(E)
 =
 2\sum_{\rho\in\mathcal R(T)}(\rho u_{n,k})(E).
$$

Now it is left to verify  the row terms. For every $\rho\in\mathcal R(T)$, the initial entry
$x_{\rho(n)}$ of the commutator part receives $e_2$.
The commutator part can be nonzero only if all its other
$\ell$ entries receive $e_1$.

When these $\ell$ entries all receive $e_1$, the commutator part equals
$$
 [e_2,e_1,\ldots,e_1]=e_2.
$$
These entries use all $\ell$ variables assigned $e_1$
under $E$, so every outer multiplier receives $e_2$.
Consequently, each row term is either zero or
$$
 L_{e_2}^{\,k}e_2\ne0.
$$
The identity row permutation gives this nonzero value.

Thus the row sum is a positive integer multiple of
$L_{e_2}^{\,k}e_2$. Since the field has characteristic
zero, this sum is nonzero. Its double is also nonzero,
which proves $f_{n,k,j}(E)\ne0$ in both cases.

Finally, the substitution $E$ assigns $e_2$ to exactly
$k+1$ variables. If $k\ge1$, every spanning $Y_1$-word
of $M_{n,k-1}$ has at most $k-1$ outer multipliers.
Its commutator part therefore receives at least two copies
of $e_2$ and vanishes. Hence $E$ forces $M_{n,k-1}$ to be zero.
For $k=0$, this is immediate from $M_{n,-1}=0$.
Since $f_{n,k,j}(E)\ne0$, we conclude that
$$
 f_{n,k,j}\notin M_{n,k-1}.
$$
\end{proof}

\begin{lemma}\label{sn-gen:layers}
For $0\le k\le n-2$, let $q=\min(k+1,n-k-1)$. Then
$$
\overline M_{n,k}
=
\bigoplus_{j=1}^{q}
\mathbb{F}S_n\bigl(f_{n,k,j}+M_{n,k-1}\bigr)
\cong
\bigoplus_{j=1}^{q}S^{(n-j,j)}.
$$

\end{lemma}

\begin{proof}

Fix $1\le j\le q$. By Lemma~\ref{sn-gen:nonvanishing},
$$
 e_{T_{n,j}}\bigl(u_{n,k}+M_{n,k-1}\bigr)
 =
 f_{n,k,j}+M_{n,k-1}
 \ne0.
$$
It follows from the Young-symmetrizer theorem that
$$
 \mathbb{F}S_n\bigl(f_{n,k,j}+M_{n,k-1}\bigr)
 \cong S^{(n-j,j)}.
$$
For distinct  $j$, these Specht modules correspond to distinct partitions and are therefore
nonisomorphic. Hence their sum is direct.

By the hook-length formula, one can have
$$
\dim S^{(n-j,j)}
=\frac{n!(n-2j+1)}{j!(n-j+1)!}
=\binom nj-\binom n{j-1}.
$$
Their direct sum consequently has dimension
$$
\sum_{j=1}^{q}
\left(\binom nj-\binom n{j-1}\right)
=\binom nq-1
=\binom n{k+1}-1.
$$
By \eqref{sn-gen:layer-dimension}, this is the dimension
of $\overline M_{n,k}$, and it proves our statement.
\end{proof}

\begin{theorem}\label{sn-gen:full-decomposition}
For every $n\ge2$,
$$
\mathcal M_n
\cong
\bigoplus_{j=1}^{\lfloor n/2\rfloor}
(n-2j+1) S^{(n-j,j)},
$$
and
$$
P_n(\mathcal V)
\cong
(n-1)S^{(n)}
\oplus
\bigoplus_{j=1}^{\lfloor n/2\rfloor}
(n-2j+1)S^{(n-j,j)}.
$$
\end{theorem}

\begin{proof}
From the filtration
$$
 0=M_{n,-1}\subseteq M_{n,0}\subseteq\cdots
 \subseteq M_{n,n-2}=\mathcal M_n,
$$ and Lemma~\ref{sn-gen:layers},
we obtain
$$
 \mathcal M_n
 \cong
 \bigoplus_{k=0}^{n-2}\overline M_{n,k}\cong \bigoplus_{k=0}^{n-2}
 \bigoplus_{j=1}^{\min(k+1,n-k-1)}
 S^{(n-j,j)}.
$$
We now collect all summands corresponding to the same
partition $(n-j,j)$.

For a fixed $j$, the summand $S^{(n-j,j)}$ appears
once for each integer $k$ satisfying
$$
 j\le k+1
 \qquad\text{and}\qquad
 j\le n-k-1.
$$
Combining them one can have
$$
 j-1\le k\le n-j-1.
$$
Since $j\ge1$ and $j\le n-j$, the possible
values of $j$ are
$$
 1\le j\le\lfloor n/2\rfloor.
$$
Consequently, interchanging the order of the two direct sums,
we obtain
$$
 \mathcal M_n
 \cong
 \bigoplus_{j=1}^{\lfloor n/2\rfloor}
 \bigoplus_{k=j-1}^{n-j-1}S^{(n-j,j)}.
$$
For each fixed $j$, the inner direct sum contains
$$
 (n-j-1)-(j-1)+1=n-2j+1
$$
copies of $S^{(n-j,j)}$. Hence
$$
 \mathcal M_n
 \cong
 \bigoplus_{j=1}^{\lfloor n/2\rfloor}
 (n-2j+1)S^{(n-j,j)}.
$$
Finally, complete reducibility and
\eqref{sn-gen:trivial-quotient} yield
$$
P_n(\mathcal V)\cong
\mathcal M_n\oplus
\bigl(P_n(\mathcal V)/\mathcal M_n\bigr)
\cong
\mathcal M_n\oplus
(n-1)S^{(n)}.
$$
\end{proof}

\section{Identities of proper subvarieties and products of commutators}
\label{comm-prod:section}\label{rect-conseq:section}
In this section we prove that every proper subvariety of $\mathcal V$ satisfies
a commutator-power identity. Then we obtain the vanishing of sufficiently long
right-normed products of commutator elements.

All powers below are right-normed: $z^1=z$, $z^{m+1}=z(z^m)$.
As before, $[a_1,\ldots,a_s]$ denotes a left-normed
Lie commutator. A polynomial described as nonzero is nonzero
modulo the identities of $\mathcal V=\operatorname{var}(B)$.

\begin{prp}\label{rect-conseq:power}
Suppose that  $\mathcal W\subseteq\mathcal V$ satisfies an identity
$h=0$, where $0\ne h\in P_n(\mathcal V)$ generates an
irreducible $S_n$-module corresponding to the two-row
Young diagram $(\lambda_1,\lambda_2)\vdash n$ with $\lambda_2>0$.
Then $\mathcal W$ satisfies $[x,y]^m=0$ for some positive
integer $m\le2\lambda_1+\lambda_2$.
\end{prp}
\begin{proof}
Set $j=\lambda_2$ and $T=T_{n,j}$.
Since $\mathbb{F}S_nh$ is irreducible of shape $(n-j,j)$,
the Young-symmetrizer theorem allows us to choose
its generator $h$ in $e_TP_n(\mathcal V)$.
This choice gives an equivalent identity. The polynomials $f_{n,k,j}$, with
$j-1\le k\le n-j-1$, form a basis of
$e_TP_n(\mathcal V)$. By 
Lemma~\ref{sn-gen:nonvanishing}, and Theorem~\ref{sn-gen:full-decomposition},
we can therefore write
\begin{equation}\label{rect-conseq:reduced-identity}
 h=\sum_{k=j-1}^{n-j-1}\alpha_k f_{n,k,j}.
\end{equation}
Assume  $K=\max\{k:\alpha_k\ne0\}$ and  set $\ell=n-K-1$. Then
$$
 h=\alpha_Kf_{n,K,j}+h_0,
 \qquad h_0\in M_{n,K-1}.
$$
Consider the evaluation on $B$ given by
$$
 E_K:\quad x_1=\cdots=x_\ell=e_1,\qquad
 x_{\ell+1}=\cdots=x_n=e_2.
$$
For $K\ge1$, every spanning word of $M_{n,K-1}$ has at most
$K-1$ outer multipliers. Its commutator part therefore
receives at least two of the $K+1$ copies of $e_2$ and
vanishes in $B$. For $K=0$, we have $M_{n,-1}=0$.
Thus $h_0(E_K)=0$.

By Lemma~\ref{sn-gen:nonvanishing} and the fact that $\alpha_K\ne0$,
there is a nonzero scalar $c\in\mathbb F$ such that
\begin{equation}\label{rect-conseq:canonical}
 h(E_K)=cL_{e_2}^{\,K}e_2.
\end{equation}
Set $z=[x,y]$ and $m=2\ell+K+1=2n-K-1$.
We will show that the following identity holds for
$\mathcal V$:
\begin{equation}\label{rect-conseq:substitution}
 h(\underbrace{z^2,\ldots,z^2}_{\ell},
   \underbrace{z,\ldots,z}_{K+1})=cz^m.
\end{equation}

Since $\mathcal V=\operatorname{var}(B)$, it suffices to
verify this equality on $B$.
Both sides are homogeneous of degree $m$ in $z$.
Since $[B,B]=\mathbb F \langle e_2\rangle$, it is enough to take $z=e_2$;
homogeneity then gives the equality for every possible
value of $z$, including zero.
For $z=e_2$, we have $z^2=e_1$, so the left-hand side is
$$
 h(E_K)=cL_{e_2}^{\,K}e_2=ce_2^{K+1}.
$$
Taking into account the multiplication table of $B$ and $m-(K+1)=2\ell$ is even, we have $e_2^{K+1}=e_2^m$.
This proves \eqref{rect-conseq:substitution} on $B$,
and consequently for $\mathcal V$.

For every $A\in\mathcal W$, substitution into the identity
$h=0$ now implies $cz^m=0$.
Since $c\ne0$, we obtain $[x,y]^m=0$.
Finally, $K\ge j-1$ implies
$$
 m=2n-K-1\le2n-j=2\lambda_1+\lambda_2.
$$
\end{proof}

By Theorem~\ref{sn-gen:full-decomposition}, the sum of the
trivial submodules of $P_n(\mathcal V)$ is isomorphic to
$(n-1)S^{(n)}$ for $n\ge2$.
Then every polynomial in this submodule is fixed by all
permutations of its variables.
Thus, a nonzero linear combination of generators of these
$n-1$ trivial modules is a nonzero symmetric multilinear
polynomial. The following proposition considers all such
linear combinations.
\begin{prp}\label{rect-conseq:trivial}
Suppose that $\mathcal W\subseteq\mathcal V$ satisfies an identity
$h=0$, where $0\ne h\in P_n(\mathcal V)$ generates an
irreducible $S_n$-module corresponding to the one-row
Young diagram $(n)$, with $n\ge2$.
Then $\mathcal W$ satisfies $[x,y]^m=0$ for some integer
$m$ with $n\le m\le2n$.
\end{prp}

\begin{proof}
Since $h\ne0$ in $P_n(\mathcal V)$ and
$\mathcal V=\operatorname{var}(B)$, the polynomial $h$
has a nonzero substitution on $B$.
By multilinearity, some substitution on basis elements is
nonzero. Since $h$ is symmetric, there is an integer
$0\le s\le n$ such that
$$
 h(\underbrace{e_1,\ldots,e_1}_{n-s},
   \underbrace{e_2,\ldots,e_2}_{s})\ne0.
$$

An induction on the degree, using the multiplication
table of $B$, shows that every bracketed product of
$e_1$ and $e_2$ belongs to $\mathbb F \langle e_1\rangle$ if the number
of occurrences of $e_2$ is even, and to $\mathbb F \langle e_2\rangle$
if it is odd.
Consequently, for some nonzero scalar $c$,
$$
 h(\underbrace{e_1,\ldots,e_1}_{n-s},
   \underbrace{e_2,\ldots,e_2}_{s})
 =cL_{e_2}^{s}e_1.$$

Again set $z=[x,y]$ and $m=2n-s$.
We will show that the following identity holds in $\mathcal V$:
\begin{equation}\label{rect-conseq:trivial-substitution}
 h(\underbrace{z^2,\ldots,z^2}_{n-s},
   \underbrace{z,\ldots,z}_{s})=cz^m.
\end{equation}
Both sides are homogeneous of degree $m$ in $z$.
As in the preceding proof, it suffices to check the equality
on $B$ with $z=e_2$.
Then $z^2=e_1$ and consequently 
$$L_{e_2}^{s}e_1=L_{e_2}^{m}e_1=L_{e_2}^{m}z^2=z^{m+2}=z^m.$$
The multiplication table of $B$
gives the required equality. Therefore it  holds for $\mathcal V$.

Since $h=0$ is an identity of $\mathcal W$,
\eqref{rect-conseq:trivial-substitution} gives $cz^m=0$.
We obtain $[x,y]^m=0$. Finally, $0\le s\le n$ implies $n\le m\le2n$.
\end{proof}

Let $A\in\mathcal V$ and set
$$
 C=[A,A]=\operatorname{span}_{\mathbb F}
 \{[a,b]:a,b\in A\}.
$$
\begin{lemma}\label{comm-prod:bicommutativity}
For all $u,v,w\in C$ and $z\in A$, we have
\begin{equation}\label{comm-prod:basic}
 uv=vu,\qquad (zv)w=0,\qquad u(vz)=v(uz).
\end{equation}
In particular, $c_1\bigl(c_2(\cdots(c_{m-1}c_m)\cdots)\bigr)$
is symmetric in $c_1,\ldots,c_m\in C$.
\end{lemma}

\begin{proof}
By the metabelian identity \eqref{metabelian},
$uv=vu$ for $u,v\in C$.
Since $z\in A$ and $v,w\in C$, by \eqref{powerful identity} we obtain
$$
 (zv)w=z[w,v]=0.
$$
Next, substitute $a=z$, $b=u$, and $c=v$ into
\eqref{mainident1}. Since $[u,v]=0$, we obtain
$$
 v(uz)-u(vz)=0.
$$
\end{proof}
 Lemma \ref{comm-prod:bicommutativity} gives the follwoing corollary.
\begin{corollary}\label{bicom}  For all $u,v,w\in C=[A,A]$, we have
$$
 (uv)w=(uw)v,\qquad u(vw)=v(uw).
$$
\end{corollary}

Recall that the following identity, obtained in the proof of Lemma~\ref{metabelian-quartic}:
\begin{equation}\label{comm-prod:plucker}
 [b,d][c,a]+[c,b][d,a]+[d,c][b,a]=0.
\end{equation}
We will use this identity in the next proposition.
\begin{prp}\label{prp1}
Let $A$ be an algebra in the variety $\mathcal V$
satisfying the identity $[x,y]^n=0$ for some integer
$n>0$. Then $A$ satisfies the identity
\begin{equation}\label{comm-prod:mixed-identity}
 [x_1,y_1]\cdots[x_n,y_n]=0,
\end{equation}
where all products are right-normed.
\end{prp}

\begin{proof}
The case $n=1$ is immediate. Suppose that $n\ge2$.
By Lemma~\ref{comm-prod:bicommutativity}, right-normed
products of elements of $[A,A]$ are symmetric in their
factors. Therefore the complete linearization of
$[x,y]^n=0$, is equivalent to
$$
 f(x_1,y_1,\ldots,x_n,y_n)
 =
 \sum_{\sigma\in S_n}
 [x_1,y_{\sigma(1)}]\cdots[x_n,y_{\sigma(n)}]=0.
$$

Set $h_0=f$ and define
$$
 h_k:=h_{k-1}(x_k,y_k)-h_{k-1}(y_k,x_k),
 \qquad 1\le k\le n,
$$
where, in the second term, $x_k$ and $y_k$ are interchanged and all other
variables are fixed. Each $h_k=0$ is an identity of $A$.

We show by induction that $h_k$ is a linear combination
of products
$$
 [x_1,y_{i_1}]\cdots[x_n,y_{i_n}],
$$
where $(i_1,\ldots,i_n)$ is a permutation of
$(1,\ldots,n)$ and $i_r=r$ for $1\le r\le k$.
For $k=0$, this follows from the definition of $f$.

Suppose that the statement holds for $h_{k-1}$, and let
$$
 z=[x_1,y_{i_1}]\cdots[x_n,y_{i_n}]
$$
be a product appearing in this linear combination.

If $i_k=k$, interchanging $x_k$ and $y_k$ reverses
the factor $[x_k,y_k]$ and leaves all other factors
unchanged. Thus alternation doubles this summand.

If $i_k\ne k$, then $z$ contains a factor $[x_m,y_k]$
for some $m>k$. By
Lemma~\ref{comm-prod:bicommutativity}, we may place
$[x_k,y_{i_k}]$ and $[x_m,y_k]$ in the last two positions
of the right-normed product. By
\eqref{comm-prod:plucker}, we have
$$
 [x_k,y_{i_k}][x_m,y_k]
 -
 [y_k,y_{i_k}][x_m,x_k]
 =
 [x_k,y_k][x_m,y_{i_k}].
$$
Applying the remaining outer multiplications preserves
this equality. Hence, alternation replaces these two
factors by $[x_k,y_k]$ and $[x_m,y_{i_k}]$, leaving
all previously fixed factors unchanged.
After reordering the factors by symmetry, the resulting
product has the required form. This completes the induction.

Finally, we obtain
$$
 h_n=\lambda[x_1,y_1]\cdots[x_n,y_n]
$$
for some nonzero $\lambda\in\mathbb F$. Since $h_n=0$ is an identity of $A$, we obtain
\eqref{comm-prod:mixed-identity}.
\end{proof} 
\begin{remark}\label{comm-prod:remark-B}
The assumption  $[x,y]^n=0$ in Proposition \eqref{prp1} cannot be omitted. In the algebra $B$, the
right-normed powers $[e_2,e_1]^n$ is $e_1$ if $n$ is even, and  $e_2$ if $n$ is odd. Thus $B$ does not satisfy $[x,y]^n=0$ for any $n$.
\end{remark}

\begin{theorem}\label{comm-prod:proper-subvarieties}
Let $\mathcal W$ be a proper subvariety of $\mathcal V$.
Then there exists an integer $s\ge1$ such that
$$
 c_1\bigl(c_2(\cdots(c_{s-1}c_s)\cdots)\bigr)=0.
$$
for every algebra $A\in\mathcal W$, where $c_1,\ldots,c_s\in [A,A]$.
The integer $s$ is the same for all $A\in\mathcal W$.
Moreover, if $\mathcal W$ has an additional nonzero
multilinear identity of degree $n\ge2$, one can choose
$s\le2n$.
\end{theorem}

\begin{proof}
Since $\mathcal W\subsetneq\mathcal V$ and the field has
characteristic zero, there is a multilinear identity
$g=0$ of $\mathcal W$ which is not an identity of
$\mathcal V$.
If its degree is one, the identity $x=0$ holds throughout
$\mathcal W$, and the conclusion holds with $s=1$.
Suppose that its degree is $n\ge2$.

Let
$$
 I_n=\{f\in P_n(\mathcal V): f=0
       \text{ is an identity of }\mathcal W\}.
$$
This is a nonzero $S_n$-submodule.
By complete reducibility, it contains an irreducible
submodule $U$. Choose $0\ne h\in U$.
Then $\mathbb{F}S_nh=U$ and $h=0$ is an identity of $\mathcal W$.

According to Theorem~\ref{sn-gen:full-decomposition}, we have
only two possibilities.
If $U\cong S^{(n)}$, the polynomial $h$ is symmetric,
and Proposition~\ref{rect-conseq:trivial} gives
$[x,y]^s=0$ for some $s\le2n$.
If $U\cong S^{(n-j,j)}$ with $j\ge1$,
Proposition~\ref{rect-conseq:power} implies the same conclusion,
with $s\le2n-j$.
Thus, in both cases, $\mathcal W$ satisfies $[x,y]^s=0$.

Finally, Proposition~\ref{prp1} implies 
$$
 [x_1,y_1]\cdots[x_s,y_s]=0.
$$
Since every element of $[A,A]$ is a linear combination of
commutators, multilinearity gives
$c_1\bigl(c_2(\cdots(c_{s-1}c_s)\cdots)\bigr)=0$ for arbitrary $c_i\in[A,A]$.
\end{proof}

\section*{Acknowledgments}

The authors are supported by grant AP26199089 from the Ministry of Science and Higher Education of the Republic of Kazakhstan.

The first author thanks Wayne State University for
its hospitality and excellent working conditions during a visit in which part of this work was carried out. 

The authors acknowledge the use of AI tools in improving the grammar of the manuscript.

\end{document}